\documentclass{article}
\usepackage{amssymb}
\usepackage{amsmath}
\usepackage{graphicx}

\newtheorem{theorem}{Theorem}

\newtheorem{lemma}{Lemma}

\newtheorem{remark}{Remark}
\numberwithin{equation}{section}
\numberwithin{lemma}{section}
\numberwithin{theorem}{section}
\numberwithin{remark}{section}
\numberwithin{corollary}{section}
\numberwithin{proposition}{section}
\numberwithin{definition}{section}
\newcommand{\dd}{{\rm d}}

\newcommand{\Oh}{{\cal O}}
\newcommand{\X}{{\rm X}}

\begin{document}

\title{On the long-time integration of stochastic gradient systems}
\author{B. Leimkuhler\footnotemark[1], C. Matthews\footnotemark[1] \ and M.V. Tretyakov\footnotemark[2]}
\maketitle

\begin{abstract}
This article addresses the weak convergence  of numerical methods for Brownian dynamics.  Typical analyses of numerical methods for stochastic differential equations focus on properties such as the weak order which estimates the asymptotic (stepsize $h\rightarrow 0$) convergence behavior of the error of finite time averages.   Recently it has been demonstrated, by study of Fokker-Planck operators, that a non-Markovian numerical method [Leimkuhler and Matthews, 2013; Leimkuhler et al., 2013] generates approximations in the long time limit with higher accuracy order (2nd order) than would be expected from its weak convergence analysis (finite-time averages are 1st order accurate).  In this article we describe the transition from the transient to the steady-state regime of this numerical method by estimating the time-dependency of the coefficients in an asymptotic expansion for the weak error, demonstrating that the convergence to 2nd order is exponentially rapid in time.  Moreover, we provide numerical tests of the theory, including comparisons of the efficiencies of the Euler-Maruyama method, the popular 2nd order Heun method, and the non-Markovian method.

\end{abstract}

\renewcommand{\thefootnote}{\arabic{footnote}} \renewcommand{\thefootnote}{%
\fnsymbol{footnote}} \footnotetext[1]{%
School of Mathematics and the Maxwell Institute for Mathematical Sciences, University of Edinburgh, Kings Buildings, Mayfield Road, Edinburgh, EH9 3JZ, UK}\footnotetext[2]{%
School of Mathematical Sciences, University of Nottingham, 
Nottingham, NG7 2RD, UK. 
Email: Michael.Tretyakov@nottingham.ac.uk}

\section{Introduction}

Stochastic gradient systems are stochastic differential equations in $d$
dimensions having the form
\begin{equation}
\mathrm{d}\mathrm{X}=a(\mathrm{X})\mathrm{d}t+\sigma \mathrm{d}\mathrm{w},\
\mathrm{X}(0)=\mathrm{X}_{0},  \label{a}
\end{equation}%
where
\begin{equation}
a(x):=-\nabla V(x),  \label{a1}
\end{equation}%
$V(x)$, $x \in {\bf R}^d$, is a potential energy function and $\sigma >0$ is a constant which
characterizes the strength of the additive noise, here described by a
standard $d$-dimensional Wiener process $\mathrm{w}(t)$. These systems
originate in the work of Einstein \cite{Ein1,Ein2} to describe the motion of
Brownian particles. They arise in mathematical models for chemistry,
physics, biology and other areas, when the cumulative effect of unresolved
degrees of freedom must be incorporated into a model to ensure its physical
relevance. Under mild conditions on $V$, the system (\ref{a}) is ergodic
\cite{Has,Stua} and has the unique invariant distribution $\rho _{\beta
}\varpropto \exp (-\beta V)$, where $\beta =2\sigma ^{-2}$. Numerical
methods for solving the equation (\ref{a}) compute a discrete sequence of
states $\mathrm{X}_{1},\mathrm{X}_{2},\ldots $ by iteratively approximating
the short time evolution. The error in the numerical solution is typically
quantified in either a \emph{strong} sense (accuracy with respect to a
particular stochastic path associated to (\ref{a})) or by reference to an
evolving distribution (\emph{weak} error, or error in averages); the latter
is the focus of this article. Ideally, the discrete states are ultimately
distributed in a way that is consistent with the invariant distribution, but
for complex applications the introduction of error in the numerical process
is inevitable. In this article we examine the asymptotic ($t\rightarrow
\infty $) behavior of the weak error.

Undoubtedly, the most common numerical method for solving (\ref{a}) is the
Euler-Maruyama method which approximates $\mathrm{X}(t_{k})$, $t_{k}=hk$, by
the iteration
\begin{equation}
\mathrm{X}_{k+1}=\mathrm{X}_{k}+ha(\mathrm{X}_{k})+\sigma \sqrt{h}\xi _{k+1},
\label{e}
\end{equation}%
where $\xi _{k}=(\xi _{k}^{1},\ldots ,\xi _{k}^{d})^{\top }$ and $\xi
_{k}^{i},$ $i=1,\ldots ,d,$ $k=1,\ldots ,$ are i.i.d. random variables with
the law $\mathcal{N}(0,1).$ For analysis of the weak error, one considers a
finite time interval $[0,\tau]$, with $\tau=hN$. The probability measure
associated to (\ref{a}) is described by a probability density $\rho (t,x)$
which evolves according to the Fokker-Planck equation
\begin{equation*}
\frac{\partial \rho }{\partial t}=\mathcal{L}^{\dagger }\rho ,
\end{equation*}%
where $\mathcal{L}^{\dagger }$ is the adjoint (in the $L_{2}$ sense) of the
generator for (\ref{a}) which is defined by
\begin{equation}
\mathcal{L}:=\sum_{i=1}^{d}a^{i}(x)\frac{\partial }{\partial x^{i}}+\frac{%
\sigma ^{2}}{2}\sum_{i=1}^{d}\frac{\partial ^{2}}{\left( \partial
x^{i}\right) ^{2}}.  \label{ee}
\end{equation}%
The solution $\rho (t,x)$ evolves from an initial probability distribution $%
\rho (0,x)$ to the steady state $\rho (\infty ,x)=\rho _{\beta }$. Let $%
\varphi $ be a test function (e.g. an element of the Schwarz space of $%
C^{\infty }$ functions rapidly decaying at infinity). Then average of $%
\varphi $ at time $\tau$ may be taken to be
\begin{equation}
\bar{\varphi}(\tau)=\mathbf{E}_{\rho (\tau,\cdot )}\varphi \equiv \int_{\mathbb{R}%
^{d}}\varphi (x)\rho (\tau,x)\mathrm{d}x.  \label{av}
\end{equation}

The discretization scheme (\ref{e}) may also be viewed as giving rise to an
evolving probability distribution, and thus one may think of the iterates in
(\ref{e}), $\mathrm{X}_{1},\mathrm{X}_{2},\ldots $, as being characterized
by densities $\rho _{1},\rho _{2},\ldots $. If stepsize $h$ is used, then
the average at time $\tau=Nh$ is given by
\begin{equation}
\hat{\varphi}(\tau,h)=\mathbf{E}_{\rho _{N}(\cdot )}\varphi \equiv \int_{%
\mathbb{R}^{d}}\varphi (x)\rho _{N}(x)\mathrm{d}x.  \label{av_d}
\end{equation}%
It is natural to compare (\ref{av}) and (\ref{av_d}) as a means of
quantifying the error as a function of $h$. We refer to this as the \emph{%
weak error}. For the Euler-Maruyama method it is known (see, e.g. \cite%
{KP,GNT04}) that
\begin{equation*}
|\bar{\varphi}(\tau)-\hat{\varphi}(\tau,h)|=\mathcal{O}(h).
\end{equation*}%
The Landau notation means that the given quantity is bounded for $%
h\rightarrow 0$ by $Ch$ where $C$ is a constant that is independent of the
stepsize. A better way to write this is
\begin{equation*}
|\bar{\varphi}(\tau)-\hat{\varphi}(\tau,h)|\leq C(\tau)h,
\end{equation*}%
since $C$ depends inherently on the time interval. This formula can be seen
as a consequence of an asymptotic expansion of the weak error, as proposed
by Talay and Tubaro \cite{TT90}. We note that $C$ also depends on the distribution of the
initial state of the system, i.e. $\rho (0,x)$, as well as the particular
observable, but we suppress these aspects in our notation. The asymptotic ($%
\tau \rightarrow \infty $) behavior of $C$ describes the performance of the
numerical method for computing averages with respect to the invariant
distribution. For the Euler-Maruyama method, one finds that $C$ is bounded
as $\tau \rightarrow \infty $, thus one obtains first order approximation of
averages both at finite time and in the long time limit.

In order to calculate averages in systems with complicated potentials and/or
a large number of variables, one often must perform numerical calculations
with a very long time interval. It is then desirable to use as large a
timestep as is reasonable in the interest of reducing the computational
effort, which is typically quantified in terms of the number of force
evaluations. Weak first-order methods like Euler-Maruyama can be inefficient
in practice. Schemes such as the second order stochastic Heun method \cite%
{KP,GNT04} can have greater efficiency: the stochastic Heun method uses two
evaluations of the force $-\nabla V$ at each timestep, thus, in comparison
to Euler-Maruyama, it must introduce less than about half the error at a
given stepsize to be deemed superior. The alternative method
discussed in this paper has been proposed in \cite{LM13}:
\begin{equation}
\mathrm{X}_{k+1}=\mathrm{X}_{k}+ha(\mathrm{X}_{k})+\sigma \frac{\sqrt{h}}{2}%
(\xi _{k}+\xi _{k+1}),  \label{b}
\end{equation}%
where $\xi _{k}=(\xi _{k}^{1},\ldots ,\xi _{k}^{i})^{\top }$ and $\xi
_{k}^{i},$ $i=1,\ldots ,d,$ $k=1,\ldots ,$ are i.i.d. random variables with
the law $\mathcal{N}(0,1)$. This method is very similar in form to the
Euler-Maruyama method (\ref{e}), and is as easy to implement, but the sums
of successive random increments are not statistically independent, so the
method is fundamentally non-Markovian in nature. The scheme was motivated in
\cite{LM13} by an analysis of Langevin dynamics algorithms. In \cite%
{LMS13}, the same method, along with some alternatives, was further analyzed
from the perspective of the invariant measure, providing a rigorous
foundation for the statement that the error in long-time averaging computed
using (\ref{b}) is of order two, i.e.
\begin{equation*}
\lim_{\tau \rightarrow \infty }|\bar{\varphi}(\tau)-\hat{\varphi}(\tau,h)|\leq
Kh^{2}.
\end{equation*}%
The remarkable feature of this estimate is that the second order accuracy is
achieved with only a single evaluation of the force at each timestep.
However, the result of \cite{LMS13} is essentially a formal analysis since
it is based entirely on the analysis of the invariant distribution and the
stationary Fokker-Planck equation. Such an operator-based approach does not
elucidate the progression from finite time averaging to infinite time
averaging and, in particular, nothing is demonstrated in \cite{LM13,LMS13}
about the weak accuracy of the method. In this article, we address this
issue, studying the way that the finite-time averages obtained using the
numerical scheme (\ref{b}) converge, as $\tau \rightarrow \infty $, to
steady-states of the numerical method. To do this, we compute the
Talay-Tubaro expansion at finite time and show that
\begin{equation*}
|\bar{\varphi}(\tau)-\hat{\varphi}(\tau,h)|\leq C_{0}(\tau)h+C_{1}(\tau)h^{2}+\ldots ,
\end{equation*}%
Then we demonstrate that
\begin{equation*}
\lim_{\tau \rightarrow \infty }C_{0}(\tau)=0,
\end{equation*}%
implying a \emph{superconvergence} property in the long-time limit.
Moreover, we show that this convergence is exponential in $\tau$.

We note that there are several recent papers (see \cite{Assyr} and
references therein), where the idea of modified differential equations is
exploited in order to construct higher-order schemes for computing ergodic limits.  This approach
provides the possibility of modifying schemes which are of weak order one on finite time
intervals to provide second order approximations in ergodic limits. However,
such modified schemes require either to evaluate derivative of forces or to
perform two force evaluations \cite{Assyr}, i.e., their computational cost
is at least as high as for the Heun scheme and substantially higher than for
(\ref{b}). Furthermore, although the theoretical approaches in our paper and in
\cite{Assyr} share some similarities, the results of \cite{Assyr} are not
applicable to the non-Markovian approximation (\ref{b}) and they do not also
include an analysis demonstrating that the leading term in the error of their modified
schemes goes to zero exponentially fast.

\section{Preliminaries}


We use the following notation for the solution of (\ref{a}): $\mathrm{X}(t)=%
\mathrm{X}_{t_{0},x}(t)$ when $\mathrm{X}(t_{0})=x,$ $t\geq t_{0},$ and also
we will write $\mathrm{X}_{x}(t)$ when $t_{0}=0.$ Recall (see, e.g. \cite%
{Has}) that the process $\mathrm{X}(t)$ is exponentially ergodic if for any $%
x\in \mathbf{R}^{d}$ and any function $\varphi $ with a polynomial growth
there are $C(x)>0$ and $\lambda >0$ such that
\begin{equation}
\left\vert \mathbf{E}\varphi (\mathrm{X}_{x}(t))-\varphi ^{erg}\right\vert
\leq C(x)e^{-\lambda t},\ \ t\geq 0,  \label{PA34}
\end{equation}%
where
\begin{equation}
\lim_{t\rightarrow \infty }\mathbf{E}\varphi (\mathrm{X}_{x}(t))=\int \varphi
(x)\rho (x)\,\mathrm{d} x:=\varphi ^{erg}.  \label{PA31}
\end{equation}%
The solution $\mathrm{X}(t)$ of (\ref{a}) is exponentially ergodic with the
Gibbs invariant density
\begin{equation*}
\rho (x)\varpropto \exp \left( -\frac{2}{\sigma ^{2}}V(x)\right)
\end{equation*}%
under the condition (see e.g. \cite{Has,Stua}): there exist $c_{0}\in
\mathbf{R}$ and $c_{1}>0$ such that
\begin{equation}
(x,a(x))\leq c_{0}-c_{1}|x|^{2}.  \label{KC2}
\end{equation}%
Under this condition, for all $p\geq 1$
\begin{equation}
\mathbf{E}|\mathrm{X}_{x}(t)|^{2p}\leq K(1+|x|^{2p}e^{-\lambda t}),
\label{momX}
\end{equation}%
where $K>0$ and $0<\lambda \leq c_{1}$ depend on $p$ (see e.g. \cite%
{Has,Stua}).

Introduce the operator $L$%
\begin{equation*}
L:=\frac{\partial }{\partial t}+\mathcal{L},
\end{equation*}%
where $\mathcal{L}$ is the generator for (\ref{a}) defined in (\ref{ee}). We
recall that the function
\begin{equation}
u(t,x)=\mathbf{E}\varphi (\mathrm{X}_{t,x}(\tau))  \label{u_def}
\end{equation}%
satisfies the Cauchy problem for the backward Kolmogorov equation
\begin{eqnarray}
Lu &=&0,  \label{u_prob} \\
u(\tau,x) &=&\varphi (x).  \notag
\end{eqnarray}%
The transition density $p(t,x,y)$ for (\ref{a}) satisfies the Fokker-Planck
(forward Kolmogorov) equation
\begin{gather}
\frac{\partial p}{\partial t}(t,x,y)=\mathcal{L}^{\dagger }p(t,x,y),\ \ t>0,
\label{PA32} \\
p(0,x,y)=\delta (y-x),  \notag
\end{gather}%
where $\mathcal{L}^{\dagger }$ is adjoint of $\mathcal{L},$ and the
invariant density $\rho (x)$ satisfies the stationary Fokker-Planck equation
\begin{equation}
\mathcal{L}^{\dagger }\,\rho (x)=0.  \label{PA33}
\end{equation}

We suppose that all components of random variables $\xi_k$ arising in (\ref%
{b}) and the Wiener process $w$ are independent. This assumption allows us
to use Ito integrals of the form $\int_{t_{k}}^{t}b(s,\mathrm{X}_{t_{k},%
\mathrm{X}_{k}}(s))\mathrm{d} \mathrm{w}(s),$ $t\geq t_{k},$ where $b(s,x)$
is a deterministic `good' function (also note that in this paper we are
considering the weak-sense convergence only). We will use the following
additional notation for this method: $\bar{\mathrm{X}}(t_{k})=\bar{\mathrm{X}%
}_{t_{k-1},\mathrm{X}_{k-1}}(t_{k})=\mathrm{X}_{k}$.

\section{Main result} \label{sec::main}

We start with a simple illustrative example.\medskip

\noindent \textbf{Example 3.1}. Let $a(x)=-\alpha x$ with $\alpha >0,$ then $%
\mathrm{X}(t)$ from (\ref{a}) is the Ornstein-Uhlenbeck process, which is
Gaussian with $\mathbf{E}\mathrm{X}_{x}(t)=xe^{-\alpha t}$ and $Cov(\mathrm{X%
}_{x}(s),\mathrm{X}_{x}(t))=\dfrac{\sigma ^{2}}{2\alpha }(e^{-\alpha
(t-s)}-e^{-\alpha (t+s)})$ for $s\leq t$. It is not difficult to calculate
that for the Euler scheme (\ref{e}):
\begin{eqnarray*}
\mathbf{E}\mathrm{X}_{N} &=&x_{0}(1-\alpha h)^{N}=x_{0}e^{-\alpha
\tau}(1+\Oh(h)),\  \\
Var(\mathrm{X}_{N}) &=&\frac{\sigma ^{2}}{2\alpha }\frac{1-(1-\alpha h)^{2N}%
}{1+\alpha h}=\frac{\sigma ^{2}}{2\alpha }(1-e^{-2\alpha \tau})-\frac{\sigma
^{2}}{2}h+e^{-2\alpha \tau}\Oh(h)+\Oh(h^{2}),\ \ \alpha h<1,
\end{eqnarray*}%
where $|\Oh(h^{p})|\leq Kh$ with $K>0$ independent of $\tau$, and for the scheme (%
\ref{b}):
\begin{eqnarray*}
\mathbf{E}\mathrm{X}_{N} &=&x_{0}(1-\alpha h)^{N}=x_{0}e^{-\alpha
\tau}(1+\Oh(h)),\  \\
Var(\mathrm{X}_{N}) &=&\frac{\sigma ^{2}}{2\alpha }\left[ 1-\frac{(1-\alpha
h)^{2N}}{1-\alpha h}\right] =\frac{\sigma ^{2}}{2\alpha }(1-e^{-2\alpha
\tau})+e^{-2\alpha \tau}\Oh(h).
\end{eqnarray*}%
We see that although both schemes have first order accuracy on finite time
intervals, the ergodic limit of the scheme (\ref{b}) is exact while the
ergodic limit of the Euler scheme approximates the ergodic limit of the
Ornstein-Uhlenbeck process with order one which is usually the case for weak
schemes of order one \cite{Tal90,MilTre07,MST10}. \medskip

In what follows we will assume the following.\medskip

\textbf{\noindent Assumption 3.1} \ \textit{The potential }$V(x)\in C^{7}(%
\mathbf{R}^{d}),$\textit{\ its first-order derivatives grow not faster than
a linear function at infinity and higher derivatives are bounded. The
relations (\ref{a1}) and (\ref{KC2}) hold. A function }$\varphi (x)\in C^{6}(%
\mathbf{R}^{d})$\textit{\ and it and its derivatives grow not faster than a
polynomial function at infinity.} \medskip

The most restrictive condition  in
Assumption~3.1 is the requirement for $a(x)=-\nabla V$ to be globally
Lipschitz:%
\begin{equation}
|a(x)|^{2}\leq K(1+|x|^{2}),  \label{globL}
\end{equation}
where $K>0$ is independent of $x\in \mathbf{R}^{d}.$
(Refer to  Remark ~\ref{rem_lip}, below, and the example presented in Subsection \ref{sec::lj_exp}.)

Introduce the multi-index $\mathbf{i}=(i_{1},\ldots ,i_{d})$ and $|\mathbf{i|%
}=\sum_{j=1}^{d}i_{j}.$ Under Assumption~3.1, we have the following. The
solution $u(t,x)$ of (\ref{u_prob}) belongs to $C^{\infty ,8}(\mathbf{R}%
_{+}\times \mathbf{R}^{d})$ and for some constant $K>0,\ \varkappa \in N,$
and $\lambda _{u}>0$ (see, e.g. \cite{Tal90})%
\begin{equation}
\left\vert u(t,x)-\varphi ^{erg}\right\vert \leq K(1+|x|^{\varkappa
})e^{-\lambda _{u}(\tau -t)},\ \ t\geq 0,  \label{up1}
\end{equation}%
and%
\begin{equation}
\left\vert \frac{\partial ^{j+|\mathbf{i|}}}{\partial ^{j}t\partial
^{i_{1}}x^{1}\cdots \partial ^{i_{d}}x^{d}}u(t,x)\right\vert \leq
K(1+|x|^{\varkappa })e^{-\lambda _{u}(\tau-t)}  \label{up2}
\end{equation}%
for all $1\leq |\mathbf{i|}\leq 8$ and $0\leq j\leq 2.$

The proof of the following lemma (which is an analogue of the moments bound (%
\ref{momX}) for the scheme (\ref{b})) is rather standard and is omitted here.

\begin{lemma}
Assume that (\ref{KC2}) and (\ref{globL}) hold. Let $\mathrm{X}_{k}$ be
defined by the scheme (\ref{b}). Then for all sufficiently small $h>0$ for
all $p\geq 1$ there is $\gamma \in (0,2c_{1})$ and $K>0$ such that
\begin{equation}
\mathbf{E}|\mathrm{X}_{k}|^{2p}\leq K(1+|x|^{2p}e^{-\gamma t_{k}}).
\label{momM}
\end{equation}
\end{lemma}

We prove the following convergence and error expansion theorem for the
scheme (\ref{b}).

\begin{theorem}
\label{thm1}Let Assumption~3.1 hold. Then the scheme (\ref{b}) is first
order weakly convergent and for all sufficiently small $h>0$ its error has
the form
\begin{equation}
\mathbf{E}\varphi (\mathrm{X}_{x}(\tau))-\mathbf{E}\varphi (\mathrm{X}%
_{N})=C_{0}(\tau ,x)h+C(\tau ,x)h^{2},  \label{thm11}
\end{equation}%
where
\begin{equation}
C_{0}(\tau,x)=\mathbf{E}\int_{0}^{\tau}B_0(t,\mathrm{X}_{x}(t))\mathrm{d} t,
\label{thm12}
\end{equation}%
\begin{eqnarray*}
B_0(t,x) &=&\frac{1}{2}\left[ \sum_{i,j=1}^{d}a^{j}(x)\frac{\partial }{%
\partial x^{j}}a^{i}(x)\frac{\partial }{\partial x^{i}}u(t,x)+\frac{\sigma
^{2}}{2}\sum_{i,j=1}^{d}\frac{\partial }{\partial x^{j}}a^{i}(x)\frac{%
\partial ^{2}}{\partial x^{i}\partial x^{j}}u(t,x)\right. \\
&&\left. +\frac{\sigma ^{2}}{2}\sum_{i,j=1}^{d}\frac{\partial ^{2}}{\left(
\partial x^{j}\right) ^{2}}a^{i}(x)\frac{\partial }{\partial x^{i}}u(t,x)%
\right] ,
\end{eqnarray*}%
and
\begin{equation*}
|C(\tau,x)|\leq K(1+|x|^{\varkappa }e^{-\lambda \tau}),
\end{equation*}%
for some $K>0,$ $\varkappa \in \mathbf{N}$ and $\lambda >0$ independent of $%
h $ and $\tau$.
\end{theorem}

\textbf{\noindent Proof.\ }Note that we shall use the letters $K,$ $%
\varkappa $ and $\lambda $ to denote various constants which are independent
of $h,$ $t,$ $\tau $, $x$. We will exploit ideas from \cite[Chapter 2]{MT} and,
in particular, from the proof of Theorem~2.2.5\ on the Talay-Tubaro
expansion. Using independence of $\mathrm{X}_{k}$ and $\mathrm{w}(t)-\mathrm{%
w}(t_{k}),$ $t\geq t_{k},$ we have
\begin{eqnarray}
R &:=&\mathbf{E}\varphi (\mathrm{X}_{x}(\tau))-\mathbf{E}\varphi (\mathrm{X}%
_{N})  \label{Db224} \\
&=&\sum_{k=0}^{N-1}\mathbf{E}(u(t_{k+1},\mathrm{X}_{t_{k},\mathrm{X}%
_{k}}(t_{k+1}))-u(t_{k+1},\bar{\mathrm{X}}_{t_{k},\mathrm{X}_{k}}(t_{k+1}))),
\notag
\end{eqnarray}%
where $u(t,x)$ is defined in (\ref{u_def}).

Expanding $u(t_{k+1},\bar{\mathrm{X}}_{t_{k},\mathrm{X}_{k}}(t_{k+1})))$ in
powers of $h$ around $\mathrm{X}_{k}$ by the usual Taylor formula, we obtain
\begin{eqnarray}
\mathbf{E} u(t_{k+1},\bar{\mathrm{X}}_{t_{k},\mathrm{X}_{k}}(t_{k+1}))) &=&%
\mathbf{E} u(t_{k+1},\mathrm{X}_{k})+\sum_{i=1}^{d}\mathbf{E}\left[ \Delta
\mathrm{X}_{k}^{i}\frac{\partial }{\partial x^{i}}u(t_{k+1},\mathrm{X}_{k})%
\right]  \label{Db226} \\
&&+\frac{1}{2}\sum_{i,j=1}^{d}\mathbf{E}\left[ \Delta \mathrm{X}%
_{k}^{i}\Delta \mathrm{X}_{k}^{j}\frac{\partial ^{2}}{\partial x^{i}\partial
x^{j}}u(t_{k+1},\mathrm{X}_{k})\right]  \notag \\
&&+\frac{1}{6}\sum_{i,j,l=1}^{d}\mathbf{E}\left[ \Delta \mathrm{X}%
_{k}^{i}\Delta \mathrm{X}_{k}^{j}\Delta \mathrm{X}_{k}^{l}\frac{\partial ^{3}%
}{\partial x^{i}\partial x^{j}\partial x^{l}}u(t_{k+1},\mathrm{X}_{k})\right]
\notag \\
&&+\frac{1}{24}\sum_{i,j,l,m=1}^{d}\mathbf{E}\left[ \Delta \mathrm{X}%
_{k}^{i}\Delta \mathrm{X}_{k}^{j}\Delta \mathrm{X}_{k}^{l}\Delta \mathrm{X}%
_{k}^{m}\frac{\partial ^{4}}{\partial x^{i}\partial x^{j}\partial
x^{l}\partial x^{m}}u(t_{k+1},\mathrm{X}_{k})\right]  \notag \\
&&+h^{3}r_{1}(t_{k},x)\,,  \notag
\end{eqnarray}%
where
\begin{equation*}
\Delta \mathrm{X}_{k}=ha(\mathrm{X}_{k})+\sigma \frac{\sqrt{h}}{2}(\xi
_{k}+\xi _{k+1})
\end{equation*}%
and
\begin{equation}
|r_{1}(t_{k},x)|\leq K(e^{-\lambda (\tau-t_{k})}+|x|^{\varkappa }e^{-\lambda \tau})
\label{thm110}
\end{equation}%
for some $K>0,$ $\varkappa \in \mathbf{N}$ and $\lambda >0$ independent of $%
h,$ $x,$ $t$ and $\tau$. To derive the estimate (\ref{thm110}), we used (\ref%
{up2}), the assumptions on $a(x)$ and its derivatives from Assumption~3.1,
and (\ref{momM}).

Introduce the auxiliary process
\begin{equation*}
\mathrm{X}_{k+1}^{\prime }=\mathrm{X}_{k}+ha(\mathrm{X}_{k})+\sigma \frac{%
\sqrt{h}}{2}\xi _{k}.
\end{equation*}%
Note that
\begin{equation*}
\mathrm{X}_{k}=\mathrm{X}_{k}^{\prime }+\sigma \frac{\sqrt{h}}{2}\xi _{k}.
\end{equation*}%
Using the Taylor expansions around $\mathrm{X}_{k}^{\prime }$, we get for
the second term in (\ref{Db226}):
\begin{gather}
\sum_{i=1}^{d}\mathbf{E}\left[ \Delta \mathrm{X}_{k}^{i}\frac{\partial }{%
\partial x^{i}}u(t_{k+1},\mathrm{X}_{k})\right] =h\sum_{i=1}^{d}\mathbf{E}
a^{i}(\mathrm{X}_{k})\frac{\partial }{\partial x^{i}}u(t_{k+1},\mathrm{X}%
_{k})+\frac{\sigma ^{2}}{4}h\sum_{i=1}^{d}\mathbf{E}\frac{\partial ^{2}}{%
\left( \partial x^{i}\right) ^{2}}u(t_{k+1},\mathrm{X}_{k}^{\prime })
\label{thm14} \\
+\frac{\sigma ^{4}}{16}h^{2}\sum_{i=1}^{d}\sum_{j=i+1}^{d}\mathbf{E}\frac{%
\partial ^{4}}{\left( \partial x^{i}\right) ^{2}\left( \partial x^{j}\right)
^{2}}u(t_{k+1},\mathrm{X}_{k}^{\prime })+\frac{\sigma ^{4}}{32}%
h^{2}\sum_{i=1}^{d}\mathbf{E}\frac{\partial ^{4}}{\left( \partial
x^{i}\right) ^{4}}u(t_{k+1},\mathrm{X}_{k}^{\prime })+h^{3}r_{2}(t_{k},x);
\notag
\end{gather}%
for the third term in (\ref{Db226}):%
\begin{gather}
\frac{1}{2}\sum_{i,j=1}^{d}\mathbf{E}\left[ \Delta \mathrm{X}_{k}^{i}\Delta
\mathrm{X}_{k}^{j}\frac{\partial ^{2}}{\partial x^{i}\partial x^{j}}%
u(t_{k+1},\mathrm{X}_{k})\right] =\frac{1}{2}h^{2}\sum_{i,j=1}^{d}\mathbf{E}%
\left[ a^{i}(\mathrm{X}_{k})a^{j}(\mathrm{X}_{k})\frac{\partial ^{2}}{%
\partial x^{i}\partial x^{j}}u(t_{k+1},\mathrm{X}_{k})\right]  \label{thm15}
\\
+\frac{\sigma ^{2}}{4}h^{2}\sum_{i,j=1}^{d}\mathbf{E}\left[ \frac{\partial }{%
\partial x^{j}}a^{i}(\mathrm{X}_{k}^{\prime })\frac{\partial ^{2}}{\partial
x^{i}\partial x^{j}}u(t_{k+1},\mathrm{X}_{k}^{\prime })\right] +\frac{\sigma
^{2}}{4}h^{2}\sum_{i,j=1}^{d}\mathbf{E}\left[ a^{i}(\mathrm{X}_{k}^{\prime })%
\frac{\partial ^{3}}{\partial x^{i}\left( \partial x^{j}\right) ^{2}}%
u(t_{k+1},\mathrm{X}_{k}^{\prime })\right]  \notag \\
+\frac{\sigma ^{2}}{4}h\sum_{i,j=1}^{d}\mathbf{E}\frac{\partial ^{2}}{\left(
\partial x^{i}\right) ^{2}}u(t_{k+1},\mathrm{X}_{k}^{\prime })+\frac{\sigma
^{4}}{8}h^{2}\sum_{i=1}^{d}\sum_{j=i+1}^{d}\mathbf{E}\frac{\partial ^{4}}{%
\left( \partial x^{i}\right) ^{2}\left( \partial x^{j}\right) ^{2}}u(t_{k+1},%
\mathrm{X}_{k}^{\prime })  \notag \\
+\frac{\sigma ^{4}}{16}h^{2}\sum_{i=1}^{d}\mathbf{E}\frac{\partial ^{4}}{%
\left( \partial x^{i}\right) ^{4}}u(t_{k+1},\mathrm{X}_{k}^{\prime
})+r_{3}(t_{k},x)h^{3};  \notag
\end{gather}%
for the fourth term in (\ref{Db226}):
\begin{gather}
\frac{1}{6}\sum_{i,j,l=1}^{d}\mathbf{E}\left[ \Delta \mathrm{X}%
_{k}^{i}\Delta \mathrm{X}_{k}^{j}\Delta \mathrm{X}_{k}^{l}\frac{\partial ^{3}%
}{\partial x^{i}\partial x^{j}\partial x^{l}}u(t_{k+1},\mathrm{X}_{k})\right]
=\frac{\sigma ^{2}}{4}h^{2}\sum_{i,j=1}^{d}\mathbf{E}\left[ a^{i}(\mathrm{X}%
_{k}^{\prime })\frac{\partial ^{3}}{\partial x^{i}\left( \partial
x^{j}\right) ^{2}}u(t_{k+1},\mathrm{X}_{k}^{\prime })\right]  \label{thm16}
\\
+\frac{\sigma ^{4}}{8}h^{2}\sum_{i=1}^{d}\sum_{j=i+1}^{d}\mathbf{E}\frac{%
\partial ^{4}}{\left( \partial x^{i}\right) ^{2}\left( \partial x^{j}\right)
^{2}}u(t_{k+1},\mathrm{X}_{k}^{\prime })+\frac{\sigma ^{4}}{16}%
h^{2}\sum_{i=1}^{d}\mathbf{E}\frac{\partial ^{4}}{\left( \partial
x^{i}\right) ^{4}}u(t_{k+1},\mathrm{X}_{k}^{\prime })+r_{4}(t_{k},x)h^{3};
\notag
\end{gather}%
for the fifth term in (\ref{Db226}):
\begin{gather}
\frac{1}{24}\sum_{i,j,l,m=1}^{d}\mathbf{E}\left[ \Delta \mathrm{X}%
_{k}^{i}\Delta \mathrm{X}_{k}^{j}\Delta \mathrm{X}_{k}^{l}\Delta \mathrm{X}%
_{k}^{m}\frac{\partial ^{4}}{\partial x^{i}\partial x^{j}\partial
x^{l}\partial x^{m}}u(t_{k+1},\mathrm{X}_{k})\right]  \label{thm17} \\
=\frac{\sigma ^{4}}{16}h^{2}\sum_{i=1}^{d}\sum_{j=i+1}^{d}\mathbf{E}\frac{%
\partial ^{4}}{\left( \partial x^{i}\right) ^{2}\left( \partial x^{j}\right)
^{2}}u(t_{k+1},\mathrm{X}_{k}^{\prime })+\frac{\sigma ^{4}}{32}%
h^{2}\sum_{i=1}^{d}\mathbf{E}\frac{\partial ^{4}}{\left( \partial
x^{i}\right) ^{4}}u(t_{k+1},\mathrm{X}_{k}^{\prime })+r_{5}(t_{k},x)h^{3}.
\notag
\end{gather}%
The functions $r_{i}(t_{k},x),$ $i=2,\ldots ,5,$ satisfy estimates of the
form (\ref{thm110}), which are derived using the same facts as in the case
of $r_{1}(t_{k},x).$

By Lemma~2.1.9 from \cite[p. 99]{MT} and again using independence of $%
\mathrm{X}_{k}$ and $\mathrm{w}(t)-\mathrm{w}(t_{k}),$ $t\geq t_{k}$, we get
\begin{equation}
\mathbf{E} u(t_{k+1},\mathrm{X}_{t_{k},\mathrm{X}_{k}}(t_{k+1}))=\mathbf{E}
u(t_{k+1},\mathrm{X}_{k})+h\mathbf{E}\mathcal{L}u(t_{k+1},\mathrm{X}_{k})+%
\frac{h^{2}}{2}\mathbf{E}\mathcal{L}^{2}u(t_{k+1},\mathrm{X}%
_{k})+r_{6}(t_{k},x)h^{3}\,.  \label{Db225}
\end{equation}%
We have for the second term in (\ref{Db225}):
\begin{eqnarray}
h\mathbf{E}\mathcal{L}u(t_{k+1},\mathrm{X}_{k}) &=&h\sum_{i=1}^{d}\mathbf{E}
a^{i}(\mathrm{X}_{k})\frac{\partial }{\partial x^{i}}u(t_{k+1},\mathrm{X}%
_{k})+\frac{\sigma ^{2}}{2}h\sum_{i=1}^{d}\frac{\partial ^{2}}{\left(
\partial x^{i}\right) ^{2}}u(t_{k+1},\mathrm{X}_{k}^{\prime })  \label{thm18}
\\
&&+\frac{\sigma ^{4}}{16}h^{2}\sum_{i,j=1}^{d}\mathbf{E}\frac{\partial ^{4}}{%
\left( \partial x^{i}\right) ^{2}\left( \partial x^{j}\right) ^{2}}u(t_{k+1},%
\mathrm{X}_{k}^{\prime })+r_{7}(t_{k},x)h^{3};  \notag
\end{eqnarray}%
for the third term in (\ref{Db225}):
\begin{gather}
\frac{h^{2}}{2}\mathbf{E}\mathcal{L}^{2}u(t_{k+1},\mathrm{X}_{k})=\frac{h^{2}%
}{2}\sum_{i,j=1}^{d}\mathbf{E} a^{i}(\mathrm{X}_{k})a^{j}(\mathrm{X}_{k})%
\frac{\partial ^{2}}{\partial x^{i}\partial x^{j}}u(t_{k+1},\mathrm{X}_{k})
\label{thm19} \\
+\frac{h^{2}}{2}\sum_{i,j=1}^{d}\mathbf{E} a^{j}(\mathrm{X}_{k})\frac{%
\partial }{\partial x^{j}}a^{i}(\mathrm{X}_{k})\frac{\partial }{\partial
x^{i}}u(t_{k+1},\mathrm{X}_{k})+\frac{\sigma ^{2}}{4}h^{2}\sum_{i,j=1}^{d}%
\mathbf{E} a^{i}(\mathrm{X}_{k}^{\prime })\frac{\partial ^{3}}{\partial
x^{i}\left( \partial x^{j}\right) ^{2}}u(t_{k+1},\mathrm{X}_{k}^{\prime })
\notag \\
+\frac{\sigma ^{2}}{4}h^{2}\sum_{i, j=1}^{d}\mathbf{E} a^{i}(\mathrm{X}%
_{k}^{\prime })\frac{\partial ^{3}}{\left( \partial x^{j}\right)
^{2}\partial x^{i}}u(t_{k+1},\mathrm{X}_{k}^{\prime })+\frac{\sigma ^{2}}{4}%
h^{2}\sum_{i,j=1}^{d}\mathbf{E}\frac{\partial ^{2}}{\left( \partial
x^{j}\right) ^{2}}a^{i}(\mathrm{X}_{k})\frac{\partial }{\partial x^{i}}%
u(t_{k+1},\mathrm{X}_{k})  \notag \\
+\frac{\sigma ^{2}}{2}h^{2}\sum_{i,j=1}^{d}\mathbf{E}\frac{\partial }{%
\partial x^{j}}a^{i}(\mathrm{X}_{k}^{\prime })\frac{\partial ^{2}}{\partial
x^{i}\partial x^{j}}u(t_{k+1},\mathrm{X}_{k}^{\prime })+\frac{\sigma ^{4}}{8}%
h^{2}\sum_{i,j=1}^{d}\mathbf{E}\frac{\partial ^{4}}{\left( \partial
x^{i}\right) ^{2}\left( \partial x^{j}\right) ^{2}}u(t_{k+1},\mathrm{X}%
_{k}^{\prime })+r_{8}(t_{k},x)h^{3}.  \notag
\end{gather}%
The functions $r_{i}(t_{k},x),$ $i=6,7,8,$ satisfy estimates of the form (%
\ref{thm110}), which are derived using the same facts as in the case of $%
r_{1}(t_{k},x)$ except $r_{6}(t_{k},x)$ where (\ref{momX}) was also used.

Let
\begin{equation*}
r(t_{k},x)=r_{6}(t_{k},x)+r_{7}(t_{k},x)+r_{8}(t_{k},x)-r_{5}(t_{k},x)-r_{4}(t_{k},x)-r_{3}(t_{k},x)-r_{2}(t_{k},x)-r_{1}(t_{k},x).
\end{equation*}%
Substituting (\ref{Db226})-(\ref{thm19}) in (\ref{Db224}), we obtain
\begin{eqnarray}
R &=&\frac{h^{2}}{2}\sum_{k=0}^{N-1}\left[ \sum_{i,j=1}^{d}\mathbf{E} a^{j}(%
\mathrm{X}_{k})\frac{\partial }{\partial x^{j}}a^{i}(\mathrm{X}_{k})\frac{%
\partial }{\partial x^{i}}u(t_{k+1},\mathrm{X}_{k})\right.  \label{finR} \\
&&\left. +\frac{\sigma ^{2}}{2}\sum_{i,j=1}^{d}\mathbf{E}\frac{\partial }{%
\partial x^{j}}a^{i}(\mathrm{X}_{k}^{\prime })\frac{\partial ^{2}}{\partial
x^{i}\partial x^{j}}u(t_{k+1},\mathrm{X}_{k}^{\prime })+\frac{\sigma ^{2}}{2}%
\sum_{i,j=1}^{d}\mathbf{E}\frac{\partial ^{2}}{\left( \partial x^{j}\right)
^{2}}a^{i}(\mathrm{X}_{k})\frac{\partial }{\partial x^{i}}u(t_{k+1},\mathrm{X%
}_{k})\right]  \notag \\
&&+\sum_{k=0}^{N-1}\mathbf{E} r(t_{k},\mathrm{X}_{k})h^{3}  \notag \\
&=&h^{2}\mathbf{E}\sum_{k=0}^{N-1}\frac{1}{2}\left[ \sum_{i,j=1}^{d}a^{j}(%
\mathrm{X}_{k})\frac{\partial }{\partial x^{j}}a^{i}(\mathrm{X}_{k})\frac{%
\partial }{\partial x^{i}}u(t_{k},\mathrm{X}_{k})+\frac{\sigma ^{2}}{2}%
\sum_{i,j=1}^{d}\frac{\partial }{\partial x^{j}}a^{i}(\mathrm{X}_{k})\frac{%
\partial ^{2}}{\partial x^{i}\partial x^{j}}u(t_{k},\mathrm{X}_{k})\right.
\notag \\
&&\left. +\frac{\sigma ^{2}}{2}\sum_{i,j=1}^{d}\frac{\partial ^{2}}{\left(
\partial x^{j}\right) ^{2}}a^{i}(\mathrm{X}_{k})\frac{\partial }{\partial
x^{i}}u(t_{k},\mathrm{X}_{k})\right] +\sum_{k=0}^{N-1}\mathbf{E} r(t_{k},%
\mathrm{X}_{k})h^{3}  \notag \\
&:=&h^{2}\mathbf{E}\sum_{k=0}^{N-1}B_0(t_{k},\mathrm{X}_{k})+\sum_{k=0}^{N-1}%
\mathbf{E} r(t_{k},\mathrm{X}_{k})h^{3},  \notag
\end{eqnarray}%
where (cf. (\ref{thm110}))%
\begin{equation}
|r(t_{k},x)|\leq K(e^{-\lambda (\tau-t_{k})}+|x|^{\varkappa }e^{-\lambda \tau}).
\label{rrr}
\end{equation}

Due to the properties of $u(t,x)$ (see ((\ref{up2}))-(\ref{up2})) and of $%
a(x)$ (see Assumption~3.1), we have
\begin{equation}
|B_0(t,x)|\leq K(1+|x|^{\varkappa })e^{-\lambda _{u}(\tau-t)}  \label{estB}
\end{equation}%
for some $K>0\ $and $\varkappa \in \mathbf{N}$ independent of $h,$ $x,$ $t$
and $\tau$. Using (\ref{estB}) and (\ref{momM}), we obtain from (\ref{finR}):
\begin{equation}
|R|\leq Kh(1+|x|^{\varkappa }e^{-\lambda \tau}),  \label{R_order}
\end{equation}%
for some constants $K>0,$ $\varkappa \in \mathbf{N}$ and $\lambda >0$
independent of $h,x,$ and $\tau$, i.e., the scheme (\ref{b}) is of first weak
order.

It remains to prove the expansion (\ref{thm11}). Consider now the $(d+1)$%
-dimensional system
\begin{eqnarray}
\mathrm{d} \mathrm{X} &=&a(\mathrm{X})\mathrm{d} t+\sigma \mathrm{d} \mathrm{%
w}(t),\;\mathrm{X}(0)=\mathrm{X}_{0\,}\,,  \label{Db229} \\
\mathrm{d} \mathrm{Y} &=&B_0(t,\mathrm{X})\mathrm{d} t,\;\mathrm{Y}(t_{0})=0\,.
\notag
\end{eqnarray}%
Solving (\ref{Db229}) by the scheme (\ref{b}), we get
\begin{equation}
\mathbf{E}\sum_{k=0}^{N-1}B_0(t_{k},\mathrm{X}_{k})h=\mathbf{E}\bar{\mathrm{Y}}%
(\tau)=\mathbf{E}\mathrm{Y}(\tau)+r_{B}(\tau,x)h=C_{0}(\tau,x)+r_{B}(\tau,x)h,
\label{Db230}
\end{equation}%
where $C_{0}(\tau,x)$ is equal to
\begin{equation}
C_{0}(\tau,x)=\mathbf{E}\mathrm{Y}(\tau)=\mathbf{E}\int_{0}^{\tau}B_0(s,\mathrm{X}%
_{x}(s))\mathrm{d} s\,  \label{Db231}
\end{equation}%
and
\begin{equation}
r_{B}(\tau,x)h=\sum_{k=0}^{N-1}\int_{t_{k}}^{t_{k+1}}\left[ \mathbf{E} B_0(s,%
\mathrm{X}_{x}(s))-\mathbf{E} B_0(t_{k},\mathrm{X}_{k})\right] \mathrm{d} s.
\label{thm111}
\end{equation}%
Introduce
\begin{equation}
\tilde{B}_0(t,x)=B_0(t,x)e^{\lambda _{u}(\tau-t)},  \label{tildaB}
\end{equation}%
for which we have (cf. (\ref{estB})):
\begin{equation*}
|\tilde{B}_0(t,x)|\leq K(1+|x|^{\varkappa }),
\end{equation*}%
where $K>0$ does not depend on $x,$ $t,$ and $\tau$. Using the demonstrated
first-order convergence of (\ref{b}) (cf. (\ref{R_order})), it is not
difficult to obtain that
\begin{eqnarray}
|r_{B}(\tau,x)|h &\leq &e^{-\lambda
_{u}(\tau-t)}\sum_{k=0}^{N-1}\int_{t_{k}}^{t_{k+1}}\left\vert \mathbf{E}\tilde{B%
}_0(s,\mathrm{X}_{x}(s))-\mathbf{E}\tilde{B}_0(t_{k},\mathrm{X}_{k})\right\vert
ds  \label{rB} \\
&\leq &e^{-\lambda _{u}(\tau-t)}h\sum_{k=0}^{N-1}\left\vert \mathbf{E}\tilde{B}_0%
(t_{k},\mathrm{X}_{x}(t_{k}))-\mathbf{E}\tilde{B}_0(t_{k},\mathrm{X}%
_{k})\right\vert +hK(1+|x|^{\varkappa }e^{-\lambda \tau})  \notag \\
&\leq &hK(1+|x|^{\varkappa }e^{-\lambda \tau}).  \notag
\end{eqnarray}%
The equality (\ref{finR}) together with (\ref{rrr}) and (\ref{Db230})-(\ref%
{rB}) implies (\ref{thm11})-(\ref{thm12}). $\ \square $\medskip

Now we prove that in the limit of $\tau \rightarrow \infty $ the scheme (\ref{b}%
) has second order of accuracy in $h.$

\begin{theorem}
\label{thm2}Let Assumption~3.1 hold. Then the coefficient $C_{0}(\tau,x)$ from (%
\ref{thm12}) goes to zero as $\tau \rightarrow \infty :\ $%
\begin{equation}
\left\vert C_{0}(\tau,x)\right\vert \leq K(1+|x|^{\varkappa })e^{-\lambda \tau}
\label{thm21}
\end{equation}%
for some constants $K>0,$ $\varkappa \in \mathbf{N}$ and $\lambda >0,$ i.e.,
over a long integration time the scheme (\ref{b}) is of order two up to
exponentially small correction.
\end{theorem}

\textbf{\noindent Proof.\ \ }We have
\begin{eqnarray}
C_{0}(\tau,x) &=&\int_{0}^{\tau}\mathbf{E} B_0(t,\mathrm{X}_{x}(t))\mathrm{d}
t=\int_{0}^{\tau}\int_{\mathbf{R}^{d}}B_0(t,y)p(t,x,y) \mathrm{d} y \mathrm{d} t
\label{thm22} \\
&=&\int_{0}^{\tau}\int_{\mathbf{R}^{d}}B_0(t,y)\rho (y)\mathrm{d} y \mathrm{d}
t+\int_{0}^{\tau}\int_{\mathbf{R}^{d}}B_0(t,y)[p(t,x,y)-\rho (y)]\mathrm{d} y
\mathrm{d} t,  \notag
\end{eqnarray}%
where $p(t,x,y)$ is the transition density for (\ref{a}) (see (\ref{PA32}))
and $\rho (y)$ is the invariant density. Using integration by parts and (\ref%
{a1}), it is not difficult to verify that for any $0\leq t\leq \tau:$
\begin{equation}
\int_{\mathbf{R}^{d}}B_0(t,y)\exp \left( -\frac{2}{\sigma ^{2}}V(y)\right)
\mathrm{d} y=0.  \label{thm23}
\end{equation}%
Further, using geometric ergodicity of $\mathrm{X}(t)$ (cf. (\ref{PA34})),
we have for $\tilde{B}_0(s,x)$ from (\ref{tildaB})%
\begin{equation}
|\mathbf{E}\tilde{B}_0(s,\mathrm{X}_{x}(t))-\int_{\mathbf{R}^{d}}\tilde{B}_0%
(s,y)\rho (y)]\mathrm{d} y|\leq K(1+|x|^{\varkappa })e^{-\lambda _{B}t},\ \
0\leq s\leq \tau,\ t>0,  \label{thm24}
\end{equation}%
for some constants $K>0,$ $\varkappa \in \mathbf{N}$ and $\lambda _{B}>0$
independent of $x,$ $t,$ and $\tau$.

Using (\ref{thm24}), we obtain for some $\lambda >0$ and all sufficiently
large $\tau >0$:
\begin{eqnarray*}
\left\vert \int_{0}^{\tau}\int_{\mathbf{R}^{d}}B_0(t,y)[p(t,x,y)-\rho
(y)]dydt\right\vert &=&\left\vert \int_{0}^{\tau}e^{-\lambda _{u}(\tau-t)}\int_{%
\mathbf{R}^{d}}\tilde{B}_0(t,y)[p(t,x,y)-\rho (y)] \mathrm{d} y \mathrm{d}
t\right\vert \\
&\leq &K(1+|x|^{\beta })e^{-\lambda \tau},
\end{eqnarray*}%
which implies (\ref{thm21}). \ $\square $

\begin{remark}
\label{rem_lip}We note that the global Lipschitz condition in Assumption~3.1
is not restrictive as the concept of rejecting exploding trajectories from
\cite{GNT04,MilTre07} can be used in implementing (\ref{b}) when the
coefficients of (\ref{a}) are not globally Lipschitz.
\end{remark}

\section{Discussion}

\textbf{1.} We emphasize that the fact that the average of $B_0(t,x)$ with respect
to the invariant measure is equal to zero (see (\ref{thm23})) is the reason
why the scheme (\ref{b}) is second order accurate in approximating
ergodic limits (see Theorem~\ref{thm2}).

\textbf{2.} In the case of the Euler scheme (\ref{e}) we get the same error
expansion as (\ref{thm11}) for the scheme (\ref{b}) but with a different $%
B_0(t,x)=B_0^{E}(t,x)$ (see \cite[Section 2.2.3]{MT}):
\begin{eqnarray*}
B_0^{E}(t,x) &=&\frac{1}{2}\left[ \sum_{i,j=1}^{d}a^{j}\frac{\partial u}{%
\partial x^{j}}a^{i}\frac{\partial u}{\partial x^{i}}+\frac{\sigma ^{2}}{2}%
\sum_{i,j}^{d}\frac{\partial ^{2}a^{j}}{\left( \partial x^{i}\right) ^{2}}%
\frac{\partial u}{\partial x^{j}}+\frac{\sigma ^{2}}{2}\sum_{i,j=1}^{d}a^{i}%
\frac{\partial ^{3}u}{\partial x^{i}\left( \partial x^{j}\right) ^{2}}\right.
\\
&&\left. +\sigma ^{2}\sum_{i,j=1}^{d}\frac{\partial a^{j}}{\partial x^{i}}%
\frac{\partial ^{2}u}{\partial x^{j}\partial x^{i}}+\frac{\sigma ^{4}}{6}%
\sum_{i,j=1}^{d}\frac{\partial ^{4}u}{\left( \partial x^{i}\right)
^{2}\left( \partial x^{j}\right) ^{2}}\right] .
\end{eqnarray*}%
The average of $B_0^{E}(t,x)$ with respect to the invariant measure is not equal
to zero and, consequently, the Euler scheme (\ref{e}) approximates ergodic
limits with order one -- the same order as its weak convergence over a
finite time interval (see also Example~3.1).

\textbf{3. }Let a one-step weak approximation $\bar{\mathrm{X}}_{t,x}(t+h)$
of the solution $\mathrm{X}_{t,x}(t+h)$ of (\ref{a}) generate a method of
order $p.$ Then, according to the Talay-Tubaro expansion \cite{TT90} (see
also \cite[Section 2.2.3]{MT}), the global error of the method has the form
\begin{equation}
R:=\mathbf{E}\varphi (\mathrm{X}_{x}(\tau))-\mathbf{E}\varphi (\bar{\mathrm{X}}%
_{x}(\tau))=C_{0}(\tau,x)h^{p}+\cdots +C_{n}(\tau,x)h^{p+n}+O(h^{p+n+1})\,,
\label{Db220}
\end{equation}%
where $n\in \mathbf{N}$ ($n$ can be arbitrarily large if the potential $V(x)$
belongs to $\mathbf{C}^{\infty }(\mathbf{R}^{d}),$ its first-order
derivatives grow not faster than a linear function at infinity\textit{\ }and
its higher derivatives of any order are bounded) and the functions $%
C_{0}(\tau,x),\ldots ,C_{n}(\tau,x)$ are independent of $h$. It follows from the
proof of Theorem~2.2.5 in \cite{MT} that the coefficients $C_{i}$ in (\ref%
{Db220}) can be presented in the form
\begin{equation*}
C_{i}(\tau,x)=\int_{0}^{\tau}\mathbf{E} B_{i}(s,\mathrm{X}_{x}(s))\mathrm{d} s.
\end{equation*}%
The function $B_{0}(s,x)$ is the coefficient at the leading term in the
one-step error expansion of the method analogous to $B(s,x)$ in Theorem~\ref%
{thm1}. The other $B_{i}(s,x),$ $i\geq 1,$ consists of the coefficient at $%
h^{p+i+1}$ from the the one-step error expansion of the method (analogously
to as $B_{0}(s,x)$ does at $h^{p+1})$ and of the coefficients at $h^{p+i+1}$
from one-step error expansions for approximations of $C_{j}$ with $j<i$ (see
details in \cite[Section 2.2.3]{MT}). Furthermore, one can deduce from the proof
of Theorem~\ref{thm2} that if the averages of $B_{i}(s,x)$ $0\leq i\leq q\leq n,$
with respect to the invariant measure are equal to zero then in the limit of
$\tau \rightarrow \infty $ the scheme has  $p+q$ order of accuracy in $h.$
Hence, such a detailed one-step error analysis is the basis for discovering
long time integration properties of numerical schemes and can serve as a
guide in the construction of highly efficient numerical methods for computing
ergodic limits for diffusions.

\section{Numerical experiments}

We compare the sampled distributions for the Euler-Maruyama scheme \eqref{e}
with the second-order (in the sense of approximating ergodic limits) scheme \eqref{b}, with both methods equal in cost
(measured in terms of evaluations of the force). We also compare the sampled
distributions with Heun's method, a second-order scheme requiring two
evaluations of $a(x) = \nabla V(x)$:
\begin{equation}  \label{h}
\begin{aligned} \hat{\X}_{k+1}&={\rm X}_{k}+ha({\rm X}_{k})+\sigma
\sqrt{h}\xi _{k+1},\\ {\rm X}_{k+1}&={\rm X}_{k}+\frac h2
\left[a(\hat{\X}_{k+1}) + a({\rm X}_{k}) \right]+\sigma \sqrt{h}\xi _{k+1}.
\end{aligned}
\end{equation}
As the scheme \eqref{b} computes exact long-time averages for all quadratic
potential energy functions $V$, it is necessary to consider anharmonic models in order to capture the representative behavior of the scheme.


\subsection{Anharmonic univariate model} \label{sec::cos_exp}

We consider solutions to \eqref{a} using the one-dimensional potential energy function
\begin{equation*}
V(x) = \cos(x),
\end{equation*}
with periodic $x \in [0,2\pi)$.

\subsubsection{Error in infinite time}

We sample the configurational distribution $\exp(-V(x))$ using trajectories generated using the Euler-Maruyama scheme \eqref{e}, Heun's method (\ref{h}) and the method \eqref{b}, where the trajectory runs over a fixed time interval of $[0,2\times10^8]$.

We note that the weak-sense convergence results are proved in Section~\ref{sec::main} under the assumption that test functions $\varphi (x)$ are sufficiently smooth and 
they and their derivatives grow not faster than polynomial functions at infinity (see Assumption~3.1). This is a usual assumption in stochastic numerics \cite{KP,MT}. At the same time, this assumption is not sufficient to guarantee convergence in distribution of the scheme \eqref{b}, which would require to consider $\varphi (x)$ being step functions. In \cite{BT} first-order weak-sense convergence of the Euler scheme and the corresponding Talay-Tubaro error expansion were proved in the case of $\varphi (x)$ being measurable bounded functions, which, in particular, implies convergence in distribution of the Euler scheme. Further, first-order convergence for density of the Euler scheme was proved in \cite{BT2}. Ideas from \cite{BT,BT2} can be exploited to extend the convergence results obtained in Section~\ref{sec::main} for the scheme \eqref{b} to include the case of nonsmooth $\varphi (x)$. Here we show and compare convergence in distribution of the scheme \eqref{b} and the other two tested methods experimentally.

\begin{figure}[ht]
\begin{center}
\includegraphics[width=\textwidth] {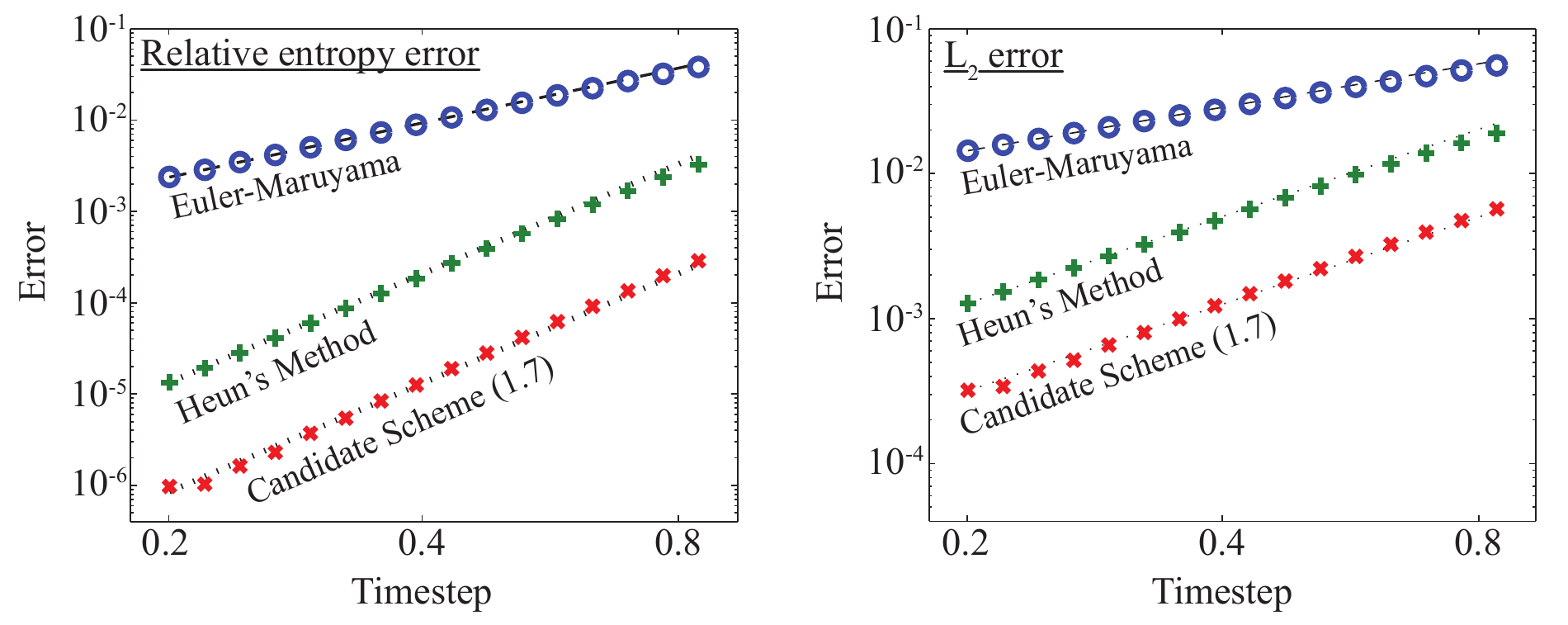}
\end{center}
\caption{The error in computed distributions is plotted for each scheme at
many stepsizes. We compare both the relative entropy (Kullback-Leibler
divergence) and the $L_2$ error of the computed distributions of $q$. The plotted black guidelines give trends with respect to stepsize, with the dashed and dotted lines giving  first and second order respectively in the right plot, and  second and fourth order respectively in the left plot.
}
\label{fig::cos1}
\end{figure}

For each scheme, we divide $[0,2\pi]$ into 100 equal histogram bins to approximate the sampled distribution, and compare the observed density of bin $i$ (denoted $\hat\rho_i$) to the exact canonical density of bin $i$ (denoted $\rho_i$) computed to high precision using a numerical solver. The error in the distribution is then reported as either the approximate $L_2$ difference in the sampled distributions, or as the relative entropy (or Kullback-Leibler divergence \cite{KL}) of the two distributions, defined by $I = \int \rho(x) \ln [ \rho(x)/\hat{\rho}(x)] \dd x$.
The relative entropy gives a measure of the information lost between two probability distributions. The two error quantities are approximated as
\begin{equation*}
 \textrm{Relative entropy error:   }\, { \sum_i \rho_i\ln\left( \frac{\rho_i}{\hat\rho_i}\right) },\qquad\qquad  \textrm{$L_2$ error:   }\,  \sqrt{ \sum_i (\hat{\rho}_i - \rho_i)^2 }. \phantom{balancing}
\end{equation*}
We compute the configurational distribution using each scheme at 16 different timesteps, where the smallest is $h=0.2$ and subsequent timesteps are increased by $10\%$. The distributions are averaged over 32 independant realizations per timestep, and the overall errors are plotted in Figure \ref{fig::cos1}.

The results match the analysis given in Section \ref{sec::main} for the large-time regime. In the case of the $L_2$ error, the Euler-Maruyama scheme gives a first order error in the computed distribution, while the other schemes give second order errors. For the computation of relative entropy, we see a doubled rate of convergence (from first to second order, or from second to fourth order).  Writing $\hat{\rho} = \rho(1+\varepsilon \psi)$, where $\varepsilon$ is a small parameter and $\int \psi \rho =0$ (conservation of total probability), we have,
\[
\int \rho\ln [1/(1+\varepsilon \psi)] \dd x = -\int \rho\ln (1+\varepsilon \psi) \dd x =
-\int \rho (\varepsilon \psi -\varepsilon^2 \psi^2 +\ldots) \dd x = -\varepsilon^2 \int \psi^2 \rho \dd x +\ldots.
\]
In the discrete context, if $\hat\rho_i = \rho_i + h^{k} \psi_i$ for an order $k$ scheme, then we find that the relative entropy is proportional to $h^{2k}$.    In practice, we observe that Heun's method and the method \eqref{b} give a fourth order relationship with the stepsize, whereas the Euler-Maruyama scheme has relative entropy proportional to $\varepsilon^2$.   The non-Markovian method gives approximately an order of magnitude improvement in this example.

\subsubsection{Error in finite time}

We consider the weak accuracy of the Euler-Maruyama scheme \eqref{e}, Heun's method (\ref{h}) and the method \eqref{b}.
In order to realize the evolving distribution computed for each scheme, we average over $2.56 \times 10^9$ independent trajectories with initial points drawn from a normal distribution with mean $\pi$ and variance 1 (where the tails of the distribution outside the periodic region are cut off). We divide $[0,2\pi]$ into 21 histogram bins, and run over $t \in [0,9]$.

\begin{figure}[ht]
\begin{center}
\includegraphics[width=5.5in] {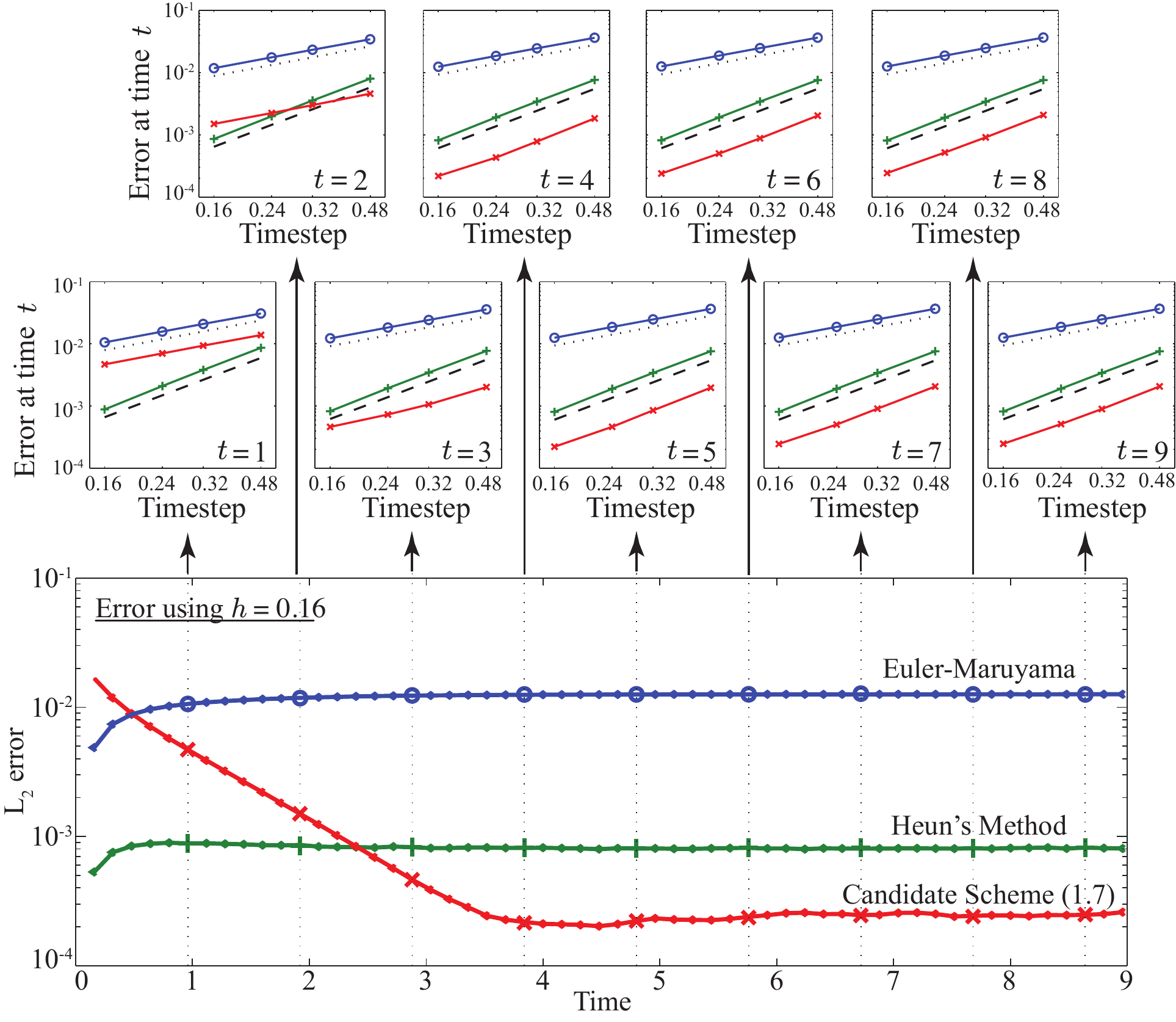}
\end{center}
\caption{The lower plot shows the error in the distribution after time $t$, as computed using each scheme at $h=0.16$. In the plots at the top, we compare the error growth with respect to stepsize $h$ at multiples of $t=0.96$. The Euler-Maruyama scheme (blue $\circ$), Heun's method (green $+$) and  the method \eqref{b} (red $\times$) are compared to first order (black, dotted) and second order (black, dashed) guidelines.}
\label{fig::cosft}
\end{figure}

As the exact solution is unknown, we compute a baseline solution using Heun's method with $h=0.04$ over the time interval. This solution is compared to the evolving distributions for $h=0.16$, $0.24$, $0.32$ and $0.48$. The growth of the error at multiples of $t=0.96$ is plotted at the top and bottom of Figure \ref{fig::cosft}, along with guidelines to indicate the order of accuracy.

We plot the error after time $t$ for each scheme, using $h=0.16$, in the central plot of Figure \ref{fig::cosft}. Initially the error in the scheme \eqref{b} reduces like $\exp(-\lambda t)$, but stabilizes after $t=4$. This is due to the behavior described in Section \ref{sec::main}, where only the first order component has an exponentially decreasing prefactor. The stabilization occurs when the $h^2$ part of the error begins to dominate the observed error.

\subsection{Lennard-Jones box} \label{sec::lj_exp}
As a more challenging problem, we compute the error in the radial distribution function for $r\in(0,6)$ for a $6\times6\times6$ periodic box of 64 Lennard-Jones particles, with interaction potential
\[
 V(q) =  \sum_{i=1}^{64}\sum_{j=i+1}^{64} r_{ij}^{-12} - r_{ij}^{-6}, \qquad r_{ij} = \|q_i - q_j\|,
\]
where $q_i$ denotes the position of particle $i$, i.e., $x$ in (\ref{a})-(\ref{a1}) is $3 \times 64 =192$-dimensional. We chose, arbitrarily,  $\beta=10$ and estimate the radial distribution function during simulation by dividing the interval $(0,6)$ into 120 histogram bins of equal length, with the error computed as the $L_2$ difference between the exact and computed radial distributions.

\begin{figure}[ht]
\begin{center}
\includegraphics[width=4.5in] {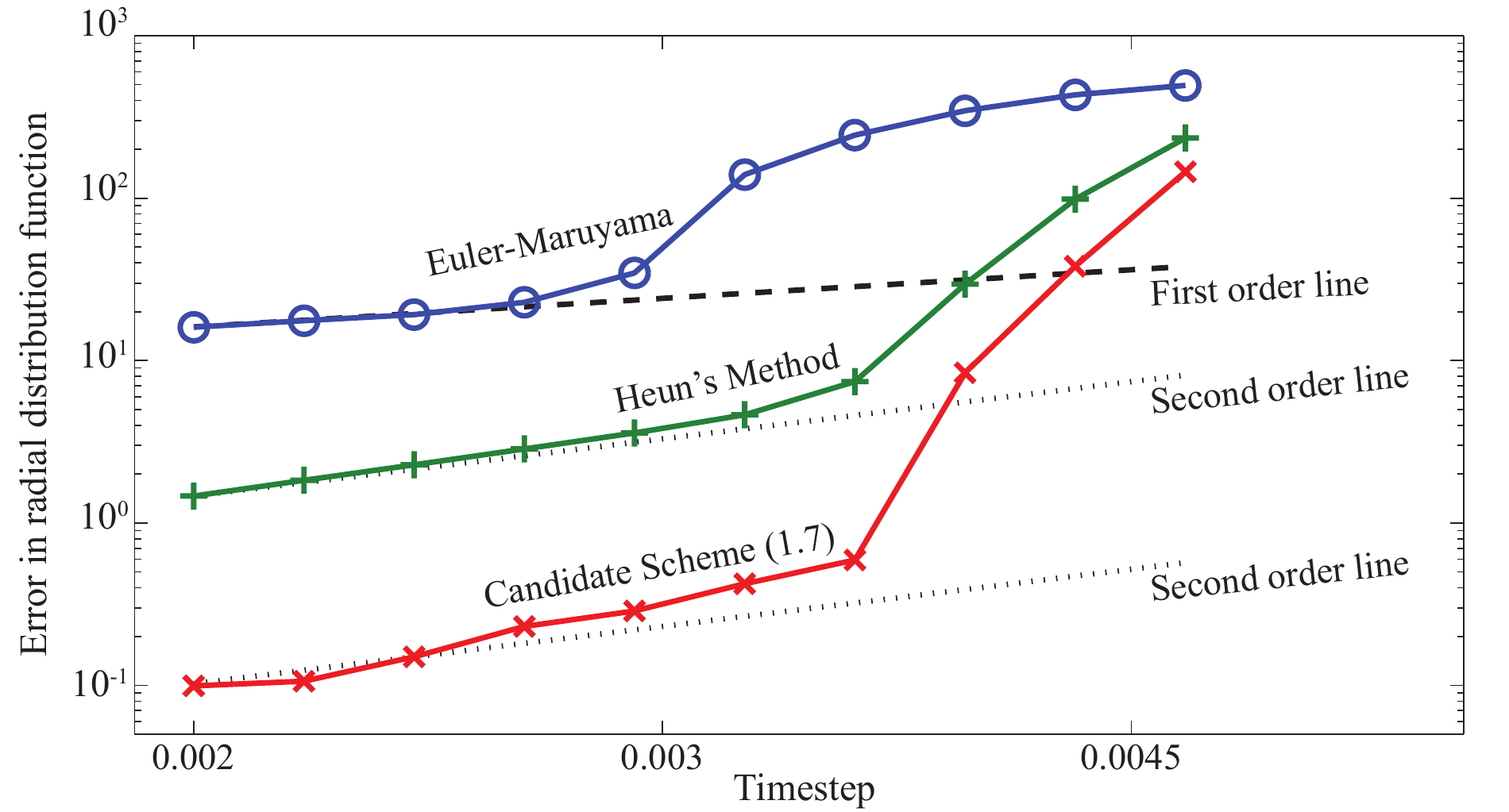}
\caption{We plot the observed $L_2$ error in the computed radial distribution functions for a periodic box of 64 Lennard-Jones particles. The Euler-Maruyama (blue $\circ$) and the method \eqref{b} (red $\times$) schemes require one force evaluation per step, while Heun's method (green $+$) requires two.} \label{fig::lj}
\end{center}
\end{figure}

We observe that the Lipschitz condition (\ref{globL}) is not, formally, satisfied for many molecular dynamics potentials (including Lennard-Jones potentials) due to the presence of singularities.   Nonetheless it is likely that, due to energetic considerations it would be possible to create a modified domain (a) in which typical solutions remain and (b) in which the Lipschitz condition (\ref{globL}) can be verified.  The numerical example presented here strongly suggests that the global Lipschitz condition could be relaxed.   More directly, the assumption (\ref{globL}) can be verified if the potential is replaced by one without singularities, e.g. by using instead Morse potentials, or by a smoothly truncated singular potential, or by a smooth Gaussian approximation of the singular potential \cite{Zi03}.

Due to the size and complexity of the problem, we cannot use standard numerical solvers to compute the exact solution. Therefore we compute a baseline solution using the scheme \eqref{b} to compute 368 realizations of a $10^7$ step trajectory (after a $10^6$ step equilibration period), with a small stepsize of $h=0.0016$.

We next compute the radial distribution functions computed using the three schemes in Section \ref{sec::cos_exp}, at ten different timesteps beginning at $h=0.002$ and with subsequent timesteps increasing by $10\%$. The trajectories were all taken over a constant time window of $[0,20000]$, with sampling beginning after a $10^6$ step equilibration.

We plot the error for all three schemes in Figure \ref{fig::lj}. For both the Euler-Maruyama scheme and Heun's method we average over 32 realizations for each timestep that we consider. This was sufficient to resolve the error introduced by these discretization methods. However, the scheme \eqref{b} proved to be sufficiently accurate that further computation was required to discern the leading error term, with the error at each timestep computed using 256 realizations to reduce the sampling error.

The results show good agreement with the theory presented in Section \ref{sec::main}. The method \eqref{b}  demonstrates an order of magnitude improvement in the long-time error of averages compared to Heun's method, while at the same time requiring half the cost (in terms of force evaluations).

\section{Summary}
In this article we have closed the gap in understanding between the typical weak error analysis of numerical discretization methods and the invariant measure accuracy of e.g. \cite{LM13,LMS13,Assyr}, demonstrating in particular that the non-Markovian numerical integration method \eqref{b} makes an exponentially rapid transition from first order weak accuracy to second order accuracy as $t\rightarrow \infty$.   Our results are confirmed in several numerical experiments, with the ultimate conclusion being that the scheme \eqref{b} is typically superior to the Euler-Maruyama and Heun's methods in terms of accuracy and efficiency for the purpose of averaging in the long term (in the transient region, the other methods may of course be better, depending on the problem).


\end{document}